\documentclass[10pt]{scrartcl}

\usepackage{amsmath, amssymb, amsthm, color, url, fullpage, supertabular}
\usepackage[colorlinks, linkcolor=black, urlcolor=blue, citecolor=black]{hyperref}
\usepackage[neverdecrease]{paralist}

\usepackage{natbib}
\setcitestyle{round}
\newcommand*{\arXiv}[1]{\bgroup\color{blue}\href{http://arxiv.org/abs/#1}{arXiv:#1}\egroup}
\newcommand*{\doi}[1]{\bgroup\color{blue}\href{http://dx.doi.org/#1}{doi:#1}\egroup}
\newcommand*{\email}[1]{\bgroup\color{blue}\href{mailto:#1}{#1}\egroup}
\renewcommand*{\url}[1]{\bgroup\color{blue}\href{#1}{#1}\egroup}
\newcommand{\todo}[1]{\bgroup\color{red}\bfseries#1\egroup}

\newcommand*{\ppara}[1]{\noindent\textbf{\textsf{#1}}\,\,}

\usepackage{mathtools}

\newcommand*{\quark}{\setbox0\hbox{$x$}\hbox to\wd0{\hss$\cdot$\hss}}

\newtheorem{theorem}{Theorem}[section]

\theoremstyle{definition}

\numberwithin{equation}{section}

\begin{document}

\title{Well-posedness of Bayesian inverse problems in quasi-Banach spaces with stable priors\footnote{To appear in the proceedings of the 88\textsuperscript{th} Annual Meeting of the International Association of Applied Mathematics and Mechanics (GAMM), Weimar 2017.  This preprint differs from the final published version in pagination and typographical detail.}}

\author{%
	T.\ J.\ Sullivan\footnote{Institute of Mathematics, Free University of Berlin, and Zuse Institute Berlin, Takustrasse 7, 14195 Berlin, Germany, \email{sullivan@zib.de}} %
}

\date{\today}

\maketitle

\begin{abstract}
	\ppara{Abstract:}
	The Bayesian perspective on inverse problems has attracted much mathematical attention in recent years.
	Particular attention has been paid to Bayesian inverse problems (BIPs) in which the parameter to be inferred lies in an infinite-dimensional space, a typical example being a scalar or tensor field coupled to some observed data via an ODE or PDE.
	This article gives an introduction to the framework of well-posed BIPs in infinite-dimensional parameter spaces, as advocated by Stuart (\emph{Acta Numer.}\ 19:451--559, 2010) and others.
	This framework has the advantage of ensuring uniformly well-posed inference problems independently of the finite-dimensional discretisation used for numerical solution.
	Recently, this framework has been extended to the case of a heavy-tailed prior measure in the family of stable distributions, such as an infinite-dimensional Cauchy distribution, for which polynomial moments are infinite or undefined.
	It is shown that analogues of the Karhunen--Lo\`eve expansion for square-integrable random variables can be used to sample such measures on quasi-Banach spaces.
	Furthermore, under weaker regularity assumptions than those used to date, the Bayesian posterior measure is shown to depend Lipschitz continuously in the Hellinger and total variation metrics upon perturbations of the misfit function and observed data.	
\end{abstract}

\section{Introduction}

The Bayesian perspective on inverse problems has attracted much mathematical attention in recent years \citep{Kaipio:2005, Stuart:2010}.
Particular attention has been paid to Bayesian inverse problems (BIPs) in which the parameter to be inferred lies in an infinite-dimensional space $\mathcal{U}$, a typical example being a scalar or tensor field coupled to some observed data via an ODE or PDE.
Numerical solution of such infinite-dimensional BIPs must necessarily be performed in an approximate manner on a finite-dimensional subspace $\mathcal{U}_{n} \subset \mathcal{U}$, but it is profitable to delay discretisation to the last possible moment and consider the original infinite-dimensional problem as the primary object of study, since infinite-dimensional well-posedness results and algorithms descend to $\mathcal{U}_{n}$ in a discretisation-independent way, whereas careless early discretisation may lead to a sequence of well-posed finite-dimensional BIPs or algorithms whose stability properties degenerate as the discretisation dimension $n$ increases.

The infinite-dimensional (or, in statistical language, \emph{non-parametric}) context presents a number of computational and theoretical challenges.
The development of sampling algorithms that are well defined in infinite-dimensional spaces, and hence robust in a high finite dimension, is an interesting topic in itself;
the preconditioned Crank--Nicolson (pCN) method of \citet{Cotter:2013} is a representative example in this area.
On a more fundamental level, though, the specification of the Bayesian prior and likelihood must be treated with care in order to ensure that the posterior measure is \emph{well posed} in the sense of Hadamard, i.e.\ is well defined and depends `nicely' upon the observed data etc.

It has become clear over the last two decades that it is not sufficient to study only the finite-dimensional formulation of BIPs, and some analysis of the infinite-dimensional limit is essential.
In the Finnish school of inverse problems, this notion of infinite-dimensional well-posedness became known as \emph{discretisation invariance}.
Its importance was highlighted by a seminal paper on edge-preserving image reconstruction using the discrete total variation prior \citep{Lassas:2004}.
The model problem in this paper was the reconstruction using $n$ pixels of a one-dimensional continuum image from $N$-dimensional linear observations subject to additive Gaussian noise.
Corresponding to a prior belief that the true image is piecewise smooth with a few jump discontinuities, the discrete total variation prior on $\mathcal{U}_{n} = \mathbb{R}^{n}$ is used.
However, \citet{Lassas:2004} showed that there is no way to consistently scale the discrete TV prior so that the Bayesian prior and posterior (and derived summary quantities such as conditional means and maximum a posteriori estimators) all have meaningful limits as $n \to \infty$.
This negative result inspired a search for priors that would provide edge-preserving and discretisation-invariant Bayesian image reconstruction.

It has recently been suggested, based upon numerical experiments, that Cauchy difference priors may be suitable priors for edge-preserving Bayesian image reconstruction \citep{Markkanen:2016}.
However, in this case, the well-posedness theory of \citet{Stuart:2010} must be extended to the heavy-tailed setting.
In this respect, this work forms part of an ongoing investigation of BIPs with heavy-tailed or other complicated structure far from the simple Gaussian regime \citep{Dashti:2012, Hosseini:2017a, Hosseini:2017b, Lassas:2009}.

\section{Bayesian inverse problems}

This section reviews the essential notions for the study of BIPs in function spaces, following the style of \citet{Stuart:2010}.
We fix two spaces $\mathcal{U}$ and $\mathcal{Y}$;
the space $\mathcal{U}$ contains an \emph{unknown} $u$ that we wish to recover from \emph{data} $y$ in $\mathcal{Y}$.
We posit that $y$ is in fact a noisily perturbed version of $G(u)$, where $G \colon \mathcal{U} \to \mathcal{Y}$ is a known function, referred to as a \emph{forward operator}.
The prototypical example to keep in mind is the case of additive Gaussian noise:
\begin{equation}
	\label{eq:forward}
	y = G(u) + \eta ,
\end{equation}
where $\eta \sim \mathcal{N}(0, \Gamma)$ is a centred $\mathcal{Y}$-valued Gaussian random variable with covariance operator $\Gamma$.
Naturally, many other noise models can arise in practice.
The challenge of recovering $u$ from $y$ is known as an \emph{inverse problem}.

Inverse problems are typically \emph{ill posed} in the sense of having no solution $u$ in the strict sense, or multiple solutions, or solutions that depend very sensitively upon the problem setup, and in particular upon the data $y$.
To deal with the existence problem, one typically relaxes the notion of solution and seeks a minimiser of a suitable misfit functional $\Phi(\setbox0\hbox{$x$}\hbox to\wd0{\hss$\cdot$\hss}; y) \colon \mathcal{U} \to \mathbb{R}$.
In the additive Gaussian case, the appropriate functional is a weighted least-squares functional:
\begin{equation}
	\label{eq:LSQ}
	\text{minimise: } u \mapsto \Phi(u; y) = \frac{1}{2} \| \Gamma^{-1/2} (G(u) - y) \|_{\mathcal{Y}}^{2} .
\end{equation}
Such problems can still be ill-conditioned, and in order to rectify this and also incorporate prior beliefs about $u$ (e.g.\ that $u$ should be small in some norm) one typically minimises a regularised version of the misfit:
\begin{equation}
	\label{eq:LSQ-R}
	\text{minimise: } u \mapsto \Phi(u; y) + R(u) .
\end{equation}
For example, on $\mathcal{U} = \mathbb{R}^{n}$, the regularisation functional $R \colon \mathcal{U} \to \mathbb{R}$ could be the Euclidean norm, corresponding to ridge regression;
or the $1$-norm, corresponding to the sparsity-promoting LASSO;
or the TV norm alluded to in the introduction.

The Bayesian interpretation of the inverse problem is to interpret both $u$ and $y$ as random variables.
One first posits a prior measure $\mu_{0}$ for $u$ --- in the case of PDE-constrained inverse problems, $\mu_{0}$ typically encodes beliefs about the smoothness of $u$ in that it charges a suitable Sobolev or H\"older space with full mass.
Eq.~\eqref{eq:forward} then determines the conditional distribution of $y|u$, the \emph{likelihood};
the Bayesian solution to the inverse problem is the conditional distribution of $u|y$, the \emph{posterior distribution}.
If $\mathcal{U}$ and $\mathcal{Y}$ were finite sets, then this posterior would be given in terms of the prior and likelihood by Bayes' rule;
informally, the posterior is the likelihood times the prior, normalised to be a probability measure.
It is also helpful to think of the prior $u$ as having probability density proportional to $\exp( -R(u))$, and the posterior $u|y$ as having density proportional to $\exp(- \Phi(u; y) - R(u))$.

However, particularly when $\mathcal{U}$ is an infinite-dimensional space, Bayes' rule must be stated more carefully, since there is no uniform reference measure (counting measure or Lebesgue measure) with respect to which one can take densities.
This generalisation of Bayes' rule is the content of Theorem~\ref{thm:mu-y-well-defined} below.
It is also of interest to study the stability of the Bayesian posterior under approximation/error in the observed data, the likelihood, or the prior;
this is the content of Theorem~\ref{thm:y-well-posed} and the discussion thereafter.

\section{Stable distributions in quasi-Banach spaces}

Cauchy distributions have recently been proposed as discretisation-invariant edge-preserving priors by \citet{Markkanen:2016}.
Cauchy distributions are an example of more general stable distributions, which are themselves generalised by infinitely divisible distributions \citep{Hosseini:2017a}.
This section surveys how stable distributions can be defined and sampled by a Karhunen--Lo\`eve-like expansion that is valid not just on Banach but even more general spaces.
While separability and completeness remain desirable attributes, it turns our that the triangle inequality can be easily relaxed.

A \emph{quasi-Banach space} is a vector space $\mathcal{U}$ equipped with a quasinorm $\| \setbox0\hbox{$x$}\hbox to\wd0{\hss$\cdot$\hss} \| \colon \mathcal{U} \to \mathbb{R}_{\geq 0}$, i.e.\ a function that is positive-definite and $1$-homogeneous and satisfies, for some choice of constant $K \geq 1$, the weakened triangle inequality
\begin{equation}
	\label{eq:weak-triangle}
	\| x + y \| \leq K ( \| x \| + \| y \| ) \text{ for all $x, y \in \mathcal{U}$,}
\end{equation}
and that is Cauchy-complete with respect to this quasinorm.
Examples of such spaces include $\ell^{p}$ and $L^{p}$ for $0 < p < 1$.

A $\mathcal{U}$-valued random variable $u$ is said to be \emph{stable of order $\alpha \in (0, 2]$} if, when $u_{1}, \dots, u_{n}$ are independent and identically distributed copies of $u$, $\sum_{i = 1}^{n} u_{i}$ is equal in distribution to $n^{1 / \alpha} u + d$ for some $d \in \mathcal{U}$.
Equivalently, in terms of the law $\mu$ of $u$ and the rescaling $\mu_{n}(E) = \mu(n^{1/\alpha} E)$,
\begin{equation}
	\label{eq:convolution}
	\mu = \underbrace{( \mu_{n} \star \dots \star \mu_{n} )}_{\text{$n$-fold convolution}} ( E + d ) 
	\text{ for all Borel-measurable $E \subseteq \mathcal{U}$.}
\end{equation}
Stability is a particularly appealing property if the aim is to construct prior measures for BIPs that are `physically consistent' in the sense of remaining in the same model class regardless of discretisation or coordinate choices, at least when the `physical quantity' obeys an additive law, e.g.\ the amount of mass or charge contained within a given physical region.

The stable distributions on $\mathbb{R}$ are completely classified by four parameters:
the index of stability $\alpha$ mentioned above, a skewness parameter $\beta \in [-1, 1]$, a scale parameter $\gamma \geq 0$, and a location parameter $\delta \in \mathbb{R}$.
We denote the unique such distribution by $\mathcal{S}(\alpha, \beta, \gamma, \delta)$.
The case $\alpha = 2$ is the case of a Gaussian (normal) distribution;
the case $\alpha = 1$, $\beta = 1$ is the Cauchy distribution.
Crucially, when $u \sim \mathcal{S}(\alpha, \beta, \gamma, \delta)$, $u$ only has finite moments up to, but not including, order $\alpha$ (the exception being $\alpha = 2$, in which case $u$ has moments of all orders).

To define and sample $\mathcal{U}$-valued stable random variables, we take a Karhunen--Lo\`eve-style approach, i.e.\ we resort to series expansions with real stable random coefficients in a countable, unconditional, normalised, Schauder basis $(\psi_{n})_{n \in \mathbb{N}}$ of $\mathcal{U}$, e.g.\ a polynomial, Fourier, or wavelet basis.
We also assume that the basis $(\psi_{n})_{n \in \mathbb{N}}$ and $q > 0$ are such that the synthesis operator $S_{\psi} \colon \underline{v} = (v_{n})_{n \in \mathbb{N}} \mapsto \sum_{n \in \mathbb{N}} v_{n} \psi_{n}$ is a continuous embedding of the sequence space $\ell^{q}$ of coefficients into $\mathcal{U}$, i.e.
\begin{equation}
	\label{eq:q-frame_upper}
	\Biggl\| \sum_{n \in \mathbb{N}} v_{n} \psi_{n} \Biggr\|_{\mathcal{U}} \leq C \| \underline{v} \|_{\ell^{q}} \text{ for all $\underline{v} \in \ell^{q}$.}
\end{equation}
When $\mathcal{U}$ is a Banach space, this assumption holds with $q = 1$ for any choice of basis $(\psi_{n})_{n \in \mathbb{N}}$, since it is just the triangle inequality for an unconditionally convergent series in $\mathcal{U}$.
Since $0 < p \leq q \leq \infty \implies \| \setbox0\hbox{$x$}\hbox to\wd0{\hss$\cdot$\hss} \|_{\ell^{q}} \leq \| \setbox0\hbox{$x$}\hbox to\wd0{\hss$\cdot$\hss} \|_{\ell^{p}}$, whenever \eqref{eq:q-frame_upper} holds for $q$ it also holds with $q$ replaced by any $p \in (0, q]$.
If inequality \eqref{eq:q-frame_upper} can be reversed, possibly with a different constant, then the basis $(\psi_{n})_{n \in \mathbb{N}}$ is known as a \emph{$q$-frame} for $\mathcal{U}$.
The case $q = 2$ is the well-known notion of a \emph{Riesz basis}.

The convergence theorem for such series expansions --- which states that sufficiently rapid decay of the scale parameters $\gamma_{n}$ implies almost sure convergence of the series in $\mathcal{U}$ --- is broadly the same as the well-known Gaussian case $\alpha = 2$, but with some logarithmic corrections \citep{Sullivan:2017}:

\begin{theorem}
	\label{thm:almost_sure_convergence}
	Consider $u = \sum_{n \in \mathbb{N}} u_{n} \psi_{n}$, where $u_{n} \sim \mathcal{S}(\alpha, \beta_{n}, \gamma_{n}, \delta_{n})$ are independent, with $\alpha \in (0, 2)$, $\underline{\beta} \subset (-1, 1)$, $\underline{\gamma} \in \ell^{\alpha}$, $\underline{\delta} \in \ell^{q}$ and, in addition,
	\begin{align}
		\label{eq:Orlicz}
		[ \underline{\gamma} ]_{\ell^{\alpha} \log \ell} & = \sum_{n \in \mathbb{N}} \bigl| \gamma_{n}^{\alpha} \log \gamma_{n} \bigr| < \infty \text{ if $\alpha = q$ or $2 q$,}
	\end{align}
	with the convention that $0 \log 0 = 0$.
	Then $u \in \mathcal{U}$ almost surely.
	Furthermore, for $0 < p \leq q$ with $p < \alpha$, $\sum_{n = 1}^{N} u_{n} \psi_{n} \to u$ in $L^{p}(\Omega, \mathbb{P}; \mathcal{U})$ as $N \to \infty$ and, in particular,
	\begin{equation}
		\label{eq:Lp_convergence_flom}
		\| u \|_{L^{p}(\Omega, \mathbb{P}; \mathcal{U})}^{p} \equiv \mathbb{E} \bigl[ \| u \|_{\mathcal{U}}^{p} \bigr] \leq C \| \underline{\gamma} \|_{\ell^{\alpha}} + C \| \underline{\delta} \|_{\ell^{q}} < \infty.
	\end{equation}
\end{theorem}

In particular, in the case $\alpha = 1$, we can define $\mathcal{U}$-valued Cauchy random variables with moments up to but not including order $1$ by using scalar Cauchy-distributed random coefficients whose scale parameters are `slightly better than $\ell^{1}$' in the sense that $\sum_{n \in \mathbb{N}} | \gamma_{n} \log \gamma_{n} |$ is finite.

\section{Well-posedness of BIPs with stable priors}

The well-posedness of the posterior $\mu^{y}$ with respect to perturbations of the problem setup is a topic of natural interest.
Since there are many choices of (inequivalent) topology or metric upon the space of probability measures on $\mathcal{U}$, this well-posedness can be quantified in many ways \citep[Chapter 14]{Deza:2016}.
We choose to use the \emph{Hellinger metric} defined by
\begin{equation}
	\label{eq:Hellinger}
	d_{\textup{H}}(\mu, \nu)
	= \frac{1}{2} \int_{\mathcal{U}} \biggl[ \sqrt{ \frac{\textup{d} \mu}{\textup{d} r} } - \sqrt{ \frac{\textup{d} \nu}{\textup{d} r} } \biggr]^{2} \, \textup{d} r
	= 1 - \mathbb{E}_{\nu} \biggl[ \sqrt{ \frac{\textup{d} \mu}{\textup{d} \nu} } \biggr] ,
\end{equation}
where $r$ is an arbitrary choice of reference measure with respect to which both $\mu$ and $\nu$ are absolutely continuous.
The Hellinger topology is stronger than the weak topology, equivalent to the total variation topology, and weaker than the relative entropy (Kullback--Leibler) topology.
It appears quite naturally in the well-posedness analysis of BIPs, and has the advantage of controlling the expected values of square-integrable quantities of interest:
\begin{equation}
	\label{eq:Hellinger-1}
	\bigl| \mathbb{E}_{\mu} [ f ] - \mathbb{E}_{\nu} [ f ] \bigr| \leq \sqrt{2} \sqrt{ \mathbb{E}_{\mu} \bigl[ |f|^{2} \bigr] + \mathbb{E}_{\nu} \bigl[ |f|^{2} \bigr] } \, d_{\textup{H}}(\mu, \nu) \text{ for all $f \in L^{2}(\mathcal{U}, \mu) \cap L^{2}(\mathcal{U}, \nu)$.}
\end{equation}

The basic result ensuring that the posterior distribution $\mu^{y}$ for $u|y$ exists and satisfies the appropriate generalisation of Bayes' rule is as follows \citep{Stuart:2010, Sullivan:2017}:

\begin{theorem}
	\label{thm:mu-y-well-defined}
	Let $\mathcal{U}$ and $\mathcal{Y}$ be quasi-Banach spaces and let $\mu_{0}$ be a prior distribution on $\mathcal{U}$.
	Suppose that $\Phi$ is locally bounded, Carath\'eodory (i.e.\ continuous in $y$ and measurable in $u$).
	Suppose also that, for every $r > 0$, there exists $M_{1. r} \colon \mathbb{R}_{\geq 0} \to \mathbb{R}$ so that
	\begin{equation}
		\label{eq:Phi-lower-bound}
		\Phi(u; y) \geq M_{1, r}(\| u \|_{\mathcal{U}}) \text{ for all $(u, y) \in \mathcal{U} \times \mathcal{Y}$ with $\| y \|_{\mathcal{Y}} < r$,}
	\end{equation}
	and such that $\exp ( - M_{1, r} ( \| \setbox0\hbox{$x$}\hbox to\wd0{\hss$\cdot$\hss} \|_{\mathcal{U}} ) ) \in L^{1} (\mathcal{U}, \mu_{0})$.
	Then, whenever $\| y \|_{\mathcal{Y}} < r$, it follows that the normalising constant $Z(y) = \mathbb{E}_{\mu_{0}} \bigl[ \exp( - \Phi( \setbox0\hbox{$x$}\hbox to\wd0{\hss$\cdot$\hss} ; y)) \bigr]$ is strictly positive and finite and
	\begin{equation}
		\frac{\textup{d} \mu^{y}}{\textup{d} \mu_{0}} (u) = \frac{\exp( - \Phi(u; y))}{Z(y)}
	\end{equation}
	does indeed define a Borel probability measure on $\mathcal{U}$, which is Radon if $\mu_{0}$ is Radon, and $\mu^{y}$ is the posterior distribution of $u \sim \mu_{0}$ conditioned upon the data $y$.
\end{theorem}

The possibility of a lower bound on $\Phi$, as in Eq.~\eqref{eq:Phi-lower-bound}, that tends to $- \infty$ as $\| u \|_{\mathcal{U}}$ or $\| y \| \to \infty$ seems counterintuitive from a finite-dimensional perspective.
However, such $\Phi$ do arise naturally in the infinite-dimensional setting, and indeed are unavoidable:
for example, in the case of a Gaussian prior and likelihood, the na{\"\i}ve formulation of $\Phi$ as a quadratic misfit would be almost surely infinite, so the Cameron--Martin theorem must be used to `subtract off the infinite part of $\Phi$', making it almost surely finite, but at the cost of a global lower bound of $- \infty$ \citep[Remark~3.8]{Stuart:2010}.

In the case of a Cauchy prior, Theorems~\ref{thm:almost_sure_convergence} and \ref{thm:mu-y-well-defined} combine to say, informally, that a BIP with Cauchy prior has well-defined posterior provided that $\Phi$ does not escape to $- \infty$ any faster than logarithmically.
This is to be contrasted with the polynomial growth rate that is permitted under a Gaussian prior.

The BIP can also be shown to be well-posed in the sense of depending continuously in the Hellinger metric upon perturbations of the data $y$ and the potential $\Phi$ \citep{Stuart:2010, Sullivan:2017}.
Indeed, with a possibly different constant, the BIP inherits the local Lipschitz continuity of $\Phi$:

\begin{theorem}
	\label{thm:y-well-posed}
	Suppose that, in addition to the conditions of Theorem~\ref{thm:mu-y-well-defined}, for every $r > 0$, there exists $M_{2, r} \colon \mathbb{R}_{\geq 0} \to \mathbb{R}$ such that, for all $(u, y, y') \in \mathcal{U} \times \mathcal{Y} \times \mathcal{Y}$ with $\| y \|_{\mathcal{Y}} , \| y ' \|_{\mathcal{Y}} < r$,
	\begin{equation}
		\label{eq:y-well-posed-in}
		\bigl| \Phi(u; y) - \Phi(u; y') \bigr| \leq \exp ( M_{2} (\| u \|_{\mathcal{U}})) \| y - y' \|_{\mathcal{Y}} .
	\end{equation}
	Suppose also that $\exp ( 2 M_{2} ( \| \setbox0\hbox{$x$}\hbox to\wd0{\hss$\cdot$\hss} \|_{\mathcal{U}} ) - M_{1, r} ( \| \setbox0\hbox{$x$}\hbox to\wd0{\hss$\cdot$\hss} \|_{\mathcal{U}} ) ) \in L^{1} (\mathcal{U}, \mu_{0})$.
	Then the posterior $\mu^{y}$ depends locally Lipschitz continuously in the Hellinger metric upon the observed data $y$ in the sense that, for all $r > 0$, there exists a constant $C > 0$ such that
	\begin{equation}
		\label{eq:y-well-posed-out}
		d_{\textup{H}} \bigl( \mu^{y}, \mu^{y'} \bigr) \leq C \| y - y' \|_{\mathcal{Y}} \text{ whenever $\| y \|_{\mathcal{Y}} , \| y ' \|_{\mathcal{Y}} < r$.}
	\end{equation}
\end{theorem}

Similar results to Theorem~\ref{thm:y-well-posed} can be shown for the well-posedness of $\mu^{y}$ with respect to perturbation of $\Phi$, e.g.\ as a result of $\Phi$ incorporating an approximate numerical solution of the forward operator $G$, which might be an ODE or PDE solution operator.
Again, similarly to \eqref{eq:y-well-posed-in}, the key property is for the approximation error in $\Phi$ to be controlled by a bound that is exponentially integrable with respect to the prior $\mu_{0}$;
the error rate for the forward problem, e.g.\ as a function of mesh size, then transfers to the BIP with a possibly different constant prefactor.

On the other hand, the well-posedness of BIPs with respect to the prior $\mu_{0}$ is a subtle topic, especially if the model is misspecified (i.e.\ there is no parameter value in $\mathcal{U}$ that corresponds to the `truth').
Of particular note in this setting is the \emph{brittleness} phenomenon highlighted by \citet{Owhadi:2015a, Owhadi:2015b}:
not only does $\mu^{y}$ depend upon $\mu_{0}$, as would be expected, but it can do so in a highly discontinuous way, in the sense that any pre-specified quantity of interest can be made to have any desired posterior expectation value after arbitrarily small perturbation in the common-moments, weak, total variation, or Hellinger topologies.
This brittleness phenomenon takes place in the limit as the data resolution becomes infinitely fine, and it is a topic of current research whether approaches such as coarsening can in general yield robust inferences \citep{Miller:2015}.

\section*{Acknowledgements}
\addcontentsline{toc}{section}{Acknowledgements}

The author is supported by the Free University of Berlin within the Excellence Initiative of the German Research FOundation (DFG).

\bibliographystyle{abbrvnat}
\addcontentsline{toc}{section}{References}
\bibliography{references.bib}

\begin{thebibliography}{14}
\providecommand{\natexlab}[1]{#1}
\providecommand{\url}[1]{\texttt{#1}}
\expandafter\ifx\csname urlstyle\endcsname\relax
  \providecommand{\doi}[1]{doi: #1}\else
  \providecommand{\doi}{doi: \begingroup \urlstyle{rm}\Url}\fi

\bibitem[Cotter et~al.(2013)Cotter, Roberts, Stuart, and White]{Cotter:2013}
S.~L. Cotter, G.~O. Roberts, A.~M. Stuart, and D.~White.
\newblock M{CMC} methods for functions: modifying old algorithms to make them
  faster.
\newblock \emph{Statist. Sci.}, 28\penalty0 (3):\penalty0 424--446, 2013.
\newblock \doi{10.1214/13-STS421}.

\bibitem[Dashti et~al.(2012)Dashti, Harris, and Stuart]{Dashti:2012}
M.~Dashti, S.~Harris, and A.~Stuart.
\newblock Besov priors for {B}ayesian inverse problems.
\newblock \emph{Inverse Probl. Imaging}, 6\penalty0 (2):\penalty0 183--200,
  2012.
\newblock \doi{10.3934/ipi.2012.6.183}.

\bibitem[Deza and Deza(2016)]{Deza:2016}
M.~M. Deza and E.~Deza.
\newblock \emph{Encyclopedia of {D}istances}.
\newblock Springer, Berlin, fourth edition, 2016.
\newblock \doi{10.1007/978-3-662-52844-0}.

\bibitem[Hosseini(2017)]{Hosseini:2017a}
B.~Hosseini.
\newblock Well-posed {B}ayesian inverse problems with infinitely divisible and
  heavy-tailed prior measures.
\newblock \emph{SIAM/ASA J. Uncertain. Quantif.}, 5\penalty0 (1):\penalty0
  1024--1060, 2017.
\newblock \doi{10.1137/16M1096372}.

\bibitem[Hosseini and Nigam(2017)]{Hosseini:2017b}
B.~Hosseini and N.~Nigam.
\newblock Well-posed {B}ayesian inverse problems: priors with exponential
  tails.
\newblock \emph{SIAM/ASA J. Uncertain. Quantif.}, 5\penalty0 (1):\penalty0
  436--465, 2017.
\newblock \doi{10.1137/16M1076824}.

\bibitem[Kaipio and Somersalo(2005)]{Kaipio:2005}
J.~Kaipio and E.~Somersalo.
\newblock \emph{Statistical and {C}omputational {I}nverse {P}roblems}, volume
  160 of \emph{Applied Mathematical Sciences}.
\newblock Springer-Verlag, New York, 2005.

\bibitem[Lassas and Siltanen(2004)]{Lassas:2004}
M.~Lassas and S.~Siltanen.
\newblock Can one use total variation prior for edge-preserving {B}ayesian
  inversion?
\newblock \emph{Inverse Problems}, 20\penalty0 (5):\penalty0 1537--1563, 2004.
\newblock \doi{10.1088/0266-5611/20/5/013}.

\bibitem[Lassas et~al.(2009)Lassas, Saksman, and Siltanen]{Lassas:2009}
M.~Lassas, E.~Saksman, and S.~Siltanen.
\newblock Discretization-invariant {B}ayesian inversion and {B}esov space
  priors.
\newblock \emph{Inverse Probl. Imaging}, 3\penalty0 (1):\penalty0 87--122,
  2009.
\newblock \doi{10.3934/ipi.2009.3.87}.

\bibitem[Markkanen et~al.(2016)Markkanen, Roininen, Huttunen, and
  Lasanen]{Markkanen:2016}
M.~Markkanen, L.~Roininen, J.~M.~J. Huttunen, and S.~Lasanen.
\newblock Cauchy difference priors for edge-preserving {B}ayesian inversion
  with an application to {X}-ray tomography, 2016.
\newblock \arXiv{1603.06135}.

\bibitem[Miller and Dunson(2015)]{Miller:2015}
J.~W. Miller and D.~B. Dunson.
\newblock Robust {B}ayesian inference via coarsening, 2015.
\newblock \arXiv{1506.06101}.

\bibitem[Owhadi et~al.(2015{\natexlab{a}})Owhadi, Scovel, and
  Sullivan]{Owhadi:2015a}
H.~Owhadi, C.~Scovel, and T.~Sullivan.
\newblock On the brittleness of {B}ayesian inference.
\newblock \emph{SIAM Rev.}, 57\penalty0 (4):\penalty0 566--582,
  2015{\natexlab{a}}.
\newblock \doi{10.1137/130938633}.

\bibitem[Owhadi et~al.(2015{\natexlab{b}})Owhadi, Scovel, and
  Sullivan]{Owhadi:2015b}
H.~Owhadi, C.~Scovel, and T.~Sullivan.
\newblock Brittleness of {B}ayesian inference under finite information in a
  continuous world.
\newblock \emph{Electron. J. Stat.}, 9\penalty0 (1):\penalty0 1--79,
  2015{\natexlab{b}}.
\newblock \doi{10.1214/15-EJS989}.

\bibitem[Stuart(2010)]{Stuart:2010}
A.~M. Stuart.
\newblock Inverse problems: a {B}ayesian perspective.
\newblock \emph{Acta Numer.}, 19:\penalty0 451--559, 2010.
\newblock \doi{10.1017/S0962492910000061}.

\bibitem[Sullivan(2017)]{Sullivan:2017}
T.~J. Sullivan.
\newblock Well-posed {B}ayesian inverse problems and heavy-tailed stable
  quasi-{B}anach space priors.
\newblock \emph{Inverse Probl. Imaging}, 11\penalty0 (5):\penalty0 857--874,
  2017.
\newblock \doi{10.3934/ipi.2017040}.

\end{thebibliography}

\end{document}